\documentclass[journal]{IEEEtran}
\usepackage{graphicx}      
\usepackage{amsmath}
\usepackage{cite}
\usepackage{multirow}
\usepackage{bm}
\usepackage{mathrsfs}
\usepackage{amssymb}
\usepackage{diagbox}
\usepackage{arydshln}
\usepackage{mathrsfs}
\usepackage{mathtools}
\usepackage[british]{babel}
\newtheorem{remark}{Remark}

\ifCLASSINFOpdf
\else

\begin{document}
%
\title{Distributed Inter-Area Oscillation Damping Control for Power Systems by Using Wind Generators and Load Aggregators}
%
%
%

\author{Zhiyuan~Tang, 
        Yue~Song, 
        Tao~Liu, 
        and~David~J.~Hill,~\IEEEmembership{Life~Fellow,~IEEE}
\thanks{This work was fully supported by the Research Grants Council of the Hong Kong Special Administrative Region under the Theme-based Research Scheme through Project No. T23-701/14-N.}
\thanks{Z. Tang, Y. Song, and T. Liu are with the Department of Electrical and Electronic Engineering, The University of Hong Kong, Hong Kong (email: zytan@eee.hku.hk; yuesong@eee.hku.hk; taoliu@eee.hku.hk).}
\thanks{D. J. Hill is with the Department of Electrical and Electronic Engineering, The University of Hong Kong, Hong Kong. He is also with the School of Electrical and Information Engineering, The University of Sydney, NSW 2006, Australia (email: dhill@eee.hku.hk; david.hill@sydney.edu.au).}}

\maketitle

\begin{abstract}
This paper investigates the potential of wind turbine generators (WTGs) and load aggregators (LAs) to provide supplementary damping control services for low frequency inter-area oscillations (LFOs) through the additional distributed damping control units (DCUs) proposed in their controllers. In order to provide a scalable methodology for the increasing number of WTGs and LAs, a novel distributed control framework is proposed to coordinate damping controllers. Firstly, a distributed algorithm is designed to reconstruct the system Jacobian matrix for each damping bus (buses with damping controllers). Thus, the critical LFO can be identified locally at each damping bus by applying eigen-analysis to the obtained system Jacobian matrix. Then, if the damping ratio of the critical LFO is less than a preset threshold, the control parameters of DCUs will be tuned in a distributed and coordinated manner to improve the damping ratio and minimize the total control cost at the same time. The proposed control framework is tested in a modified IEEE 39-bus test system. The simulation results with and without the proposed control framework are compared to demonstrate the effectiveness of the proposed framework.
\end{abstract}


\IEEEpeerreviewmaketitle

\section{Introduction}

Low frequency inter-area oscillations (LFOs) have always been a matter of concern to power system operators due to their potential threats to the power system stability \cite{kundur}. With the development of the electricity market and growing power demand, future power systems will become more stressed and operate closer to their stability limits, which highlights the need to improve the damping ratio of LFOs and prevent sustained oscillations that can result in serious consequences such as system separations or even large-area blackouts \cite{kundur}.

The power system stabilizers (PSSs) installed on conventional synchronous machines are the most important components to improve system damping against LFOs. However, the increasing penetration of wind power limits the availability of PSSs to provide sufficient damping against LFOs. For one thing, the displacement of conventional synchronous generators with wind turbine generators (WTGs) may reduce the damping ratio of inter-area modes by the reconfiguration of line power flows, reduction of system inertia, and interaction of converter controls with power system dynamics \cite{Tie_line}. For another thing, once the conventional synchronous machines are replaced by WTGs, the associated PSSs are also removed from the system with no replacement controllers for WTGs to provide damping control services. Thus, if no new alternative controllers are developed to provide supplementary damping control services, insufficient system controls may jeopardize the system security and stability. To solve this issue, in this paper, we are looking for solutions from both the generation and load sides.

For the generation side, we utilize the converter interfaced WTGs which can provide damping torques for LFOs by quickly adjusting their active power outputs though a proper control of electronic devices that interface them with the grid \cite{WPP,M_WPPs}. For the load side, the option of using highly distributed controllable loads (demand control) is appealing. Due to properties such as instantaneous responses and spatial distributions, demand control has gained a lot of attention \cite{Nondisruptive,Hierarchical,angle}. In particular, demand control has been utilized to accomplish important system support tasks such as frequency control \cite{Nondisruptive}, voltage control \cite{Hierarchical}, and small-disturbance angle stability enhancement \cite{angle}. However, the ability of demand control to provide supplementary damping control services agaist LFOs has not been thoroughly investigated yet. In this paper, the load aggregators (LAs) will be coordinated with WTGs to provide damping torques against LFOs through the additional distributed damping control units (DCUs) developed in their controllers.

In the literature, numerous methods have been proposed to coordinate traditional damping controllers (e.g. PSS) \cite{reviews,LMI2,H2} and new damping controllers (e.g. FACTS and HVDC) \cite{BMI,MPC,fuzzy1}. Approaches based on robust control theories and linear matrix inequalities have been utilized to deal with the uncertainties of operating conditions \cite{LMI2,H2,BMI}. For example, in \cite{H2}, the synthesis of the controller is formulated as a mixed $H_2$/$H_\infty$ output feedback control problem with regional pole placement that is resolved through a linear matrix inequality approach. However, such a robust controller design method is too conservative and unable to incorporate all system constraints (e.g. hard limits on the control signals). Approaches based on model predictive control have been utilized to incorporate all system constraints. For example, the authors of \cite{MPC} propose a model predictive control based HVDC supplementary controller which can incorporate plant constraints explicitly. Unfortunately, the model used in such a method is developed at a pre-given operating point, and hence, the obtained controller cannot directly guarantee robustness around the other operating points. Approaches based on fuzzy logic have been utilized to handle the variations of operating points \cite{fuzzy1}. For example, a fuzzy logic adaptive control unit is proposed in \cite{fuzzy1} to adjust control gains for different operating points. However, this fuzzy logic based method becomes very complicated when the number of damping controllers becomes large. Moreover, all the methods mentioned above are carried out in a centralized manner that lacks scalability and flexibility, i.e., a new damping controller is added into the original control system, the whole control law need to be redesigned.

To overcome the drawbacks of the abovementioned methods, in this paper, a novel distributed control framework is proposed to coordinate damping controllers, which can be implemented by local measurements and limited communications between neighboring buses. The proposed distributed control framework consists of two modules: a critical LFO identification module and a controller parameters tuning module where the communication network used in each module is different. The critical LFO identification module aims at reconstructing the system Jacobian matrix for each damping bus (a bus with damping controller) in a distributed manner where the communication network used covers all buses in the system. Thus, the critical LFO (the LFO with the least damping ratio) can be identified locally at each damping bus by applying eigen-analysis to the obtained system Jacobian matrix. Further, if the damping ratio of the critical LFO is less than a preset threshold, the parameters of DCUs will be tuned in a distributed manner to improve the damping ratio of the critical LFO and minimize the total control cost at the same time where the communication network used only covers those damping buses. The contributions of this paper are listed below:
\begin{itemize}
\item A novel two-step communication based distributed control framework is proposed to coordinate LAs and WTGs. The proposed control method can survive one-point failure in the communication network and is suitable in practice for its scalability.
\item In the critical LFO identification module, based on structural properties of the original power grid, a distributed calculation algorithm is developed to recover the Jacobian matrices for each damping bus.
\item In the controller parameters tuning module, based on the eigenvalue sensitivities, a controller tuning problem is formulated and solved in a distributed manner.
\end{itemize}

The rest of the paper is organized as follows. Section \uppercase\expandafter{\romannumeral2} introduces the DCU and the power system model to be studied. The proposed distributed control framework is explicitly presented in Section \uppercase\expandafter{\romannumeral3}. Section \uppercase\expandafter{\romannumeral4} presents a case study by using a modified IEEE 39-bus test system. Conclusions are given in Section \uppercase\expandafter{\romannumeral5}.

\subsection*{Notations}

Denote $\mathbb{R}$ and $\mathbb{C}$ as the set of real numbers and complex numbers, respectively. An $m$-dimensional vector is denoted as $\bm x=[x_i]\in\mathbb{R}^m$. The transpose of a vector or a matrix is defined as $(\cdot)^T$. The notation $\bm I_m\in \mathbb{R}^{m\times m}$ denotes the identity matrix, $\bm 0$ is a zero vector or matrix with an appropriate dimention, and $\bm e_i\in \mathbb{R}^p$ denotes the vector with the $i^{th}$ entry being one and others being zeros. The notation $|x|$ ($\angle x$) takes the modulus (angle) of a complex number $x\in\mathbb{C}$. The notation $\mathbb{V}(\bm A)$ means converting the matrix $\bm A=[\bm a_1,\ldots,\bm a_p]\in \mathbb{R}^{m\times p}$ with $\bm a_i\in\mathbb{R}^m$ with $i=1,\ldots,p$ to a vector, i.e., $\mathbb{V}(\bm A)=[\bm a_1^T,\ldots,\bm a_p^T]^T\in \mathbb{R}^{mp}$. The symbols $\|\cdot\|$ and $\|\cdot\|_\infty$ denote the $l_2$ and $l_\infty$ norms for a vector, respectively.

\section{Network Description}

In this section, we firstly introduce the DCU proposed for each damping controller. Then, the power system network to be studied is introduced, which will be used to design the control framework in Section \uppercase\expandafter{\romannumeral3}.

\subsection{Distributed Damping Control Unit}

\begin{figure}[!t]
\begin{center}
\includegraphics[width=2.8in,height=0.6in]{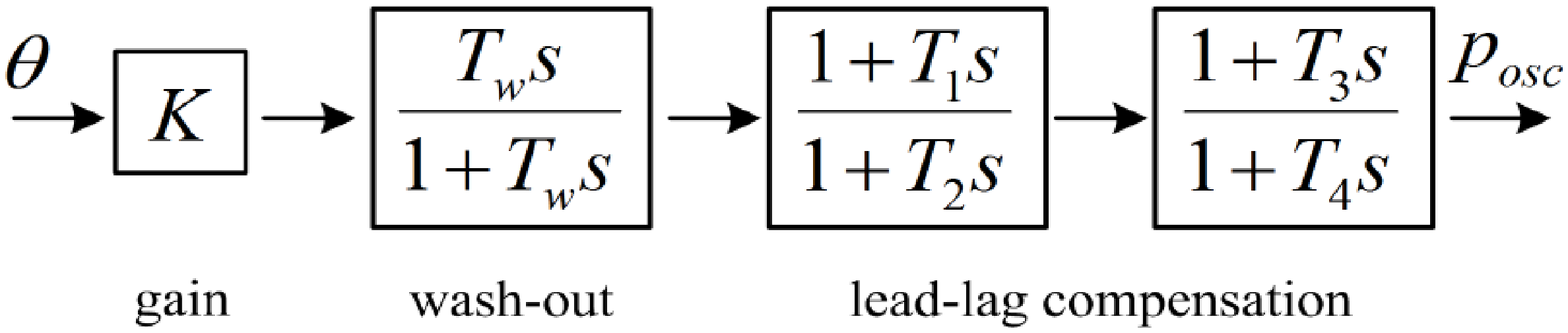}    
\caption{The control block diagram of the proposed DCU.}
\label{fig:DDCU}
\end{center}
\end{figure}

Fig. \ref{fig:DDCU} shows the block diagram of the proposed DCU which mimics the structure of PSS. The input is the local bus voltage angle $\theta_i$, and the output is $p_{osci}$ which is added to the reference active power demand of the WTG or LA to provide supplementary damping control services. The mathematical model of the $i^{th}$ DCU is given by
\begin{equation}\label{eq:DDCU}
\begin{split}
&\dot{x}_{1i}=-\frac{1}{T_{wi}}\left(K_i\theta_i+x_{1i}\right)\\
&\dot{x}_{2i}=\frac{1}{T_{2i}}\left((1-\frac{T_{1i}}{T_{2i}})(K_i\theta_i+x_{1i})-x_{2i}\right)\\
&\dot{x}_{3i}=\frac{1}{T_{4i}}\left((1-\frac{T_{3i}}{T_{4i}})\left(x_{2i}+\left(\frac{T_{1i}}{T_{2i}}(K_i\theta_i+x_{1i})\right)\right)-x_{3i}\right)\\
&p_{osci}=x_{3i}+\frac{T_{3i}}{T_{4i}}\left(x_{2i}+\frac{T_{1i}}{T_{2i}}(K_i\theta_i+x_{1i})\right).
\end{split}
\end{equation}
The dynamics can be written in a compact form as $\dot{\bm x}_{Ci}=\bm f_{Ci}(\bm x_{Ci},\theta_i)$ where $\bm x_{Ci}=[x_{1i},x_{2i},x_{3i}]^T$ is the supplementary state variables, $K_i$ is the gain, $T_{wi}$ is the wash-out time constant, $T_{1i}$, $T_{2i}$, $T_{3i}$, and $T_{4i}$ are time constants for lead-lag compensation. In the proposed control framework, $K_i$, $T_{1i}$, $T_{2i}$, $T_{3i}$, and $T_{4i}$ will be tuned to improve the damping ratio of the critical LFO.

\subsection{Power system network}

Consider a connected power system consisting of $N$ buses with $N_G$ synchronous generators (SGs), $N_W$ WTGs, $N_L$ loads, and $N_T$ transfer buses where $N=N_G+N_W+N_L+N_T$. The SG (WTG or load) bus refers to a bus that connects a SG (WTG or load) only. The transfer bus is a bus with no generation or demand. We number the SG buses as $\mathcal{V}_G=\{1,\ldots,N_G\}$, WTG buses as $\mathcal{V}_W=\{N_G+1,\ldots,N_G+N_W\}$, load buses as $\mathcal{V}_L=\{N_G+N_W+1,\ldots,N_G+N_W+N_L\}$, and transfer buses as $\mathcal{V}_T=\{N_G+N_W+N_L+1,\ldots,N\}$.

\subsubsection{SG model}

To highlight the effectiveness of the proposed damping controllers, PSSs are not included in the SG models. With the $4^{th}$-order two-axis synchronous machine model and IEEE standard exciter model (IEEET1), the mathematical model of the $i^{th}$ SG is written as:
\begin{equation}\label{eq:si}
\begin{split}
\dot{\bm x}_{Gi}&=\bm f_{Gi}(\bm x_{Gi},\theta_i,v_i)\\
p_{Gi}&=g_{p_{Gi}}(\bm x_{Gi},\theta_i,v_i)\\
q_{Gi}&=g_{q_{Gi}}(\bm x_{Gi},\theta_i,v_i),~i\in\mathcal{V}_G
\end{split}
\end{equation}
where the state variable $\bm x_{Gi}$ is defined as $\bm x_{Gi}=[e'_{qi},e'_{di},\delta_i,\omega_i,x_{mi},x_{r1i},x_{r2i},x_{fi}]^T$; $e'_{qi}$ and $e'_{di}$ are transient d-axis and q-axis voltages, respectively; $\delta_i$ and $\omega_i$ are the rotor angle and speed, respectively; $x_{mi}$, $x_{r1i}$, $x_{r2i}$ and $x_{fi}$ are the state variables corresponding to the IEEET1 exciter. The algebraic variables are the local bus voltage angle $\theta_i$ and magnitude $v_i$. The active and reactive power injections of the $i^{th}$ SG bus are denoted as $p_{Gi}$ and $q_{Gi}$, respectively. The detailed descriptions of nonlinear functions $\bm f_{Gi}$, $g_{p_{Gi}}$, $g_{q_{Gi}}$ can be found in \cite{stability_control}, which is given in the Appendix A for self-completeness.

\subsubsection{WTG model}

\begin{figure}[!t]
\begin{center}
\includegraphics[width=2.6in,height=1.0in]{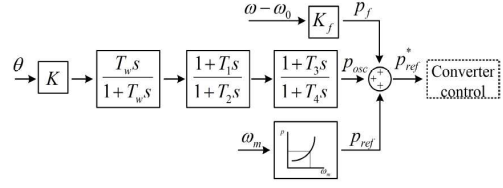}    
\caption{Control block diagram of the controller for WTG.}
\label{fig:WTG}
\end{center}
\end{figure}

Fully rated converter WTGs are adopted, which employ the configuration of a synchronous machine with a permanent magnet rotor \cite{wind_turbine}. Normally, the controller of WTG gives a reference active power demand to optimize the wind energy capture based on the measured rotor speed (see the lower branch in Fig. \ref{fig:WTG}). In this paper, two additional control units are added into the original WTG's controller to adapt the active power reference set point, i.e., the primary frequency support unit proposed in \cite {wind_frequency} (see the upper branch in Fig. \ref{fig:WTG}) and the DCU (see the middle branch in Fig. \ref{fig:WTG}). The mathematical model of the $i^{th}$ WTG is written as:
\begin{equation} \label{eq:ci}
\begin{split}
\dot{\theta}_i&=\omega_i\\
\dot{\bm x}_{Wi}&=\bm f_{Wi}(\bm x_{Wi},\omega_i,\theta_i,v_i)\\
p_{Wi}&=g_{p_{Wi}}(\bm x_{Wi},\omega_i,\theta_i,v_i)\\
q_{Wi}&=g_{q_{Wi}}(\bm x_{Wi},\omega_i,\theta_i,v_i),~i\in\mathcal{V}_W
\end{split}
\end{equation}
where $\omega_i$ is the local bus frequency. The state variable $\bm x_{Wi}=[\omega_{mi},\theta_{pi},i_{sqi},i_{cdi},\bm x_{Ci}^T]^T$ where $\omega_{m_i}$ is the rotor speed; $\theta_{p_i}$ is the pitch angle used for maximum power control; $i_{qsi}$ is the generator stator quadrature current used for active power/speed control; and $i_{cdi}$ is the converter direct current used for reactive power/voltage control; $\bm x_{Ci}=[x_{1i},x_{2i},x_{3i}]^T$ are state variables corresponding to the DCU. The active and reactive power injections of the $i^{th}$ WTG bus are denoted as $p_{Wi}$ and $q_{Wi}$, respectively. The detailed descriptions of nonlinear functions $\bm f_{Wi}$, $g_{p_{Wi}}$, $g_{q_{Wi}}$ are given in the Appendix B.

\subsubsection{Load model}

The active power of each load $p_{Li}$ is divided into two parts: the controllable part $d_i=p_{osci}(\bm x_{Li},\theta_i)$ (referred to as LA in this paper) and static voltage frequency dependent part, whereas the reactive power of each load $q_{Li}$ is assumed to be static voltage frequency dependent. With the additional DCU, the mathematical model of the $i^{th}$ load bus is given as follows:
\begin{equation} \label{eq:Ll}
\begin{split}
\dot{\theta}_i=&\omega_i\\
\dot{\bm x}_{Li}=&\bm f_{Li}(\bm x_{Li},\theta_i)\\
p_{Li}=&p_{oi}(v_i)^{\alpha_i}(1+k_{pfi}(\omega_i-\omega_0))+d_i(\bm x_{Li},\theta_i)\\
\coloneqq&g_{p_{Li}}(\bm x_{Li},\omega_i,\theta_i,v_i)\\
q_{Li}=& q_{o_i}(v_i)^{\beta_i}(1+k_{qf_i}(\omega_i-\omega_0))\\
\coloneqq&g_{q_{Li}}(\omega_i,v_i),~i\in\mathcal{V}_L
\end{split}
\end{equation}
where the state variable $\bm x_{Li}=[x_{1i},x_{2i},x_{3i}]^T$ corresponds to the DCU; $p_{oi}$ and $q_{oi}$ are the nominal values; $\alpha_i$ and $\beta_i$ are voltage coefficients; $k_{pfi}$ and $k_{qfi}$ are frequency coefficients.

\subsubsection{Transfer bus}

As transfer buses have no generations or loads, the $i^{th}$ transfer bus is simply modeled as:
\begin{equation} \label{eq:Tl}
p_{Ti}=0,~q_{Ti}=0,~i\in\mathcal{V}_T
\end{equation}
where $p_{Ti}$ and $q_{Ti}$ are the active and reactive power injections, respectively.

\subsubsection{Network power flows}

The network power flows are represented by the usual set of algebraic power flow equations, which are used to couple all buses power injection equations mentioned above. For the $i^{th}$ bus in the system, the power flow equations are given as:
\begin{equation} \label{eq:gpq}
\begin{split}
0&=-p_{i}^{inj}+v_i\sum_{j=1}^Nv_j(G_{ij}\cos\theta_{ij}+B_{ij}\sin\theta_{ij})\\
0&=-q_{i}^{inj}+v_i\sum_{j=1}^Nv_j(G_{ij}\sin\theta_{ij}-B_{ij}\cos\theta_{ij}),~i\in\mathcal{V}
\end{split}
\end{equation}
where $G_{ij}$ and $B_{ij}$ are the real and imaginary parts of $Y_{ij}$ which is the $(i,j)$ entry of the admittance matrix $\bm Y$; the notation $\theta_{ij}$ is the short for $\theta_i-\theta_j$; the set $\mathcal{V}=\mathcal{V}_G\cup\mathcal{V}_W\cup\mathcal{V}_L\cup\mathcal{V}_T$; $p_{i}^{inj}$ and $q_{i}^{inj}$ are injected active and reactive power of the $i^{th}$ bus, respectively. In particular, for SG buses, $p_{i}^{inj}=p_{Gi}$ and $q_{i}^{inj}=q_{Gi}$; for WTG buses, $p_{i}^{inj}=p_{Wi}$ and $q_{i}^{inj}=q_{Wi}$; for load buses, $p_{i}^{inj}=-p_{Li}$ and $q_{i}^{inj}=-q_{Li}$; and for transfer buses, $p_{i}^{inj}=p_{Ti}$ and $q_{i}^{inj}=q_{Ti}$.
\subsubsection{Overall system}

Combining (\ref{eq:si})-(\ref{eq:gpq}), the overall system can be expressed as differential-algebraic equations:
\begin{equation} \label{eq:fh}
\begin{split}
\dot{\bm x}&=\bm f(\bm x, \bm y)\\
\bm 0&=\bm h(\bm x, \bm y)
\end{split}
\end{equation}
where the vector $\bm x=[\bm x_G^T,\bm \theta_W^T,\bm x_W^T,\bm \theta_L^T,\bm x_L^T]^T$ and the vector $\bm y=[\bm \theta_G^T,\bm \omega_W^T,\bm \omega_L^T,\bm \theta_T^T,\bm v_G^T,\bm v_W^T,\bm v_L^T,\bm v_T^T]^T$; $\bm x_k=[\bm x_{k1}^T,\ldots,\bm x_{ki}^T,\ldots,\bm x_{k{N_k}}^T]^T,~i\in \mathcal{V}_k$, $k\in\{G,W,L\}$; $\bm \theta_k=[\theta_i]\in\mathbb{R}^{N_k},~i\in \mathcal{V}_k$, $k\in\{G,W,L,T\}$; $\bm v_k=[v_i]\in\mathbb{R}^{N_k},~i\in \mathcal{V}_k$, $k\in\{G,W,L,T\}$; $\bm \omega_k=[\omega_i]\in\mathbb{R}^{N_k},~i\in \mathcal{V}_k$, $k\in\{W,L\}$. The nonlinear functions $\bm f$ and $\bm h$ represent the system dynamics and network power flow equations, respectively.


Linearizing system (\ref{eq:fh}) gives the following linear model:
\begin{equation} \label{eq:ss}
\left [\begin{array}{c}
\Delta \dot{\bm x}\\
\hdashline[2pt/2pt]
\bm 0
\end{array} \right]
=
\left [\begin{array}{c;{2pt/2pt}c}
\bm A_{s} & \bm B_{s}\\
\hdashline[2pt/2pt]
\bm C_{s} & \bm D_{s}\\
\end{array} \right]
\left [\begin{array}{c}
\Delta \bm x\\
\hdashline[2pt/2pt]
\Delta \bm y
\end{array} \right]
\end{equation}
where the detailed expressions of the matrices $\bm A_s$, $\bm B_s$, $\bm C_s$, and $\bm D_s$ are given in the Appendix C.
Assuming $\bm D_s$ is nonsingular (it is a common assumption adopted in the literature \cite{fully_micro}) and eliminating $\Delta\bm y$ in (\ref{eq:ss}) gives:
\begin{equation} \label{eq:re}
\Delta\dot{\bm x}=\bm A_r \Delta\bm x
\end{equation}
where $\bm A_r=\bm A_s-\bm B_s\bm D_s^{-1}\bm C_s\in \mathbb{R}^{N_A\times N_A}$ with $N_A=8N_G+8N_W+4N_L$.


\section{Distributed Control Framework}

In this section, the critical LFO identification module and controller parameters tuning module that form the proposed distributed control framework will be introduced in details.

\subsection{Critical LFO identification module}

As mentioned earlier, this module aims at identifying the critical LFO for each damping bus in a distributed manner. It is known that the critical LFO can be investigated by applying eigenvalue analysis based on the global information, i.e., the system Jacobian matrix $\bm A_r$ in (\ref{eq:re}) which is usually obtained in a centralized manner \cite{kundur}. However, in this paper, we will reconstruct the matrix $\bm A_r$ for each damping bus in a distributed manner by revealing the structure properties of the power grid contained in the matrices $\bm A_s$, $\bm B_s$, $\bm C_s$, and $\bm D_s$.

By performing an elementary column operation, the matrices $\bm A_s$, $\bm B_s$, $\bm C_s$, and $\bm D_s$ can be reformulated as:
\begin{equation} \label{eq:reform}
\left [\begin{array}{c;{2pt/2pt}c}
\bm A_{s} & \bm B_{s}\\
\hdashline[2pt/2pt]
\bm C_{s} & \bm D_{s}
\end{array} \right]
=
\left [\begin{array}{c;{2pt/2pt}c}
\bm K_1 & \bm K_2\\
\hdashline[2pt/2pt]
\bm K_3 & \bm J_{pf}
\end{array} \right]\bm T
\end{equation}
where the matrix $\bm T$ is the elementary column operator and $\bm J_{pf}$ is the power flow Jacobian matrix. The detailed expressions of the matrices $\bm T$, $\bm K_1$, $\bm K_2$, $\bm K_3$, and $\bm J_{pf}$ are given in the Appendix C.

Through the matrix transform (\ref{eq:reform}), we can see that the matrices $\bm A_s$, $\bm B_s$, $\bm C_s$, and $\bm D_s$ can be reconstructed by all Jacobian matrices $\bm K_\vee^\wedge$ (refer to (\ref{eq:TMatrix}) in the Appendix C for details), $\bm J_{pf}$, identity matrices, and $\bm T$. For identity matrices and $\bm T$, since they are constant, they can be broadcasted or stored at each damping bus in advance. For all Jacobian matrices $\bm K_\vee^\wedge$ (all are block diagonal matrices) and $\bm J_{pf}$, we adopt the  distributed algorithm proposed in \cite{fully_micro} that has total $2N$ steps to calculate their elements, where the communication network used covers all buses in the system and has the same topology as the physical grid. The communication network can be described by the undirected graph $\mathcal{G}_1=\{\mathcal{V},\mathcal{E}\}$, where $\mathcal{V}$ is the set of nodes (buses) and $\mathcal{E}\subseteq\mathcal{V}\times\mathcal{V}$ represents the set of edges (branches). The set of neighbors of node $i$ is represent by $\mathcal{N}_i=\{j\in \mathcal{V}:(j,i)\in \mathcal{E}\}$ with cardinality $|\mathcal{N}_i|=\mathcal{D}_i$. We assume that 1) each bus knows the parameters of its local machine (or load) and lines connecting it; 2) each damping bus knows the model structure of SG, WTG, and load; 3) each bus knows its own bus number and total number of buses $N$; 4) each bus in the system has the capability of local measurement, storing data, processing data, communicating with its neighbors, and calculation; and 5) communication delays are negligible.

At each step, bus $i,~i\in \mathcal{V}$ has four columns of data for communication, denoted as $\bm \gamma_i^a$, $\bm \varpi_i^a$, $\bm \gamma_i^b$, $\bm \varpi_i^b\in \mathbb{R}^{2N}$. The data update process is designed as follows:
\begin{equation} \label{eq:update}
[\bm X^a(\tau+1),\bm X^b(\tau+1)]=\bm J_{pf} [\bm X^a(\tau),\bm X^b(\tau)]
\end{equation}
where $\bm X^a(\tau)$, $\bm X^b(\tau)\in \mathbb{R}^{2N\times 2N}$ are the data matrices at the $\tau^{th}$ step iteration with the definitions as follows:
\begin{equation} \label{eq:X}
\begin{split}
\bm X^a(\tau)&=[\bm \gamma_1^a(\tau),\ldots,\bm \gamma_{N}^a(\tau),\bm \varpi_1^a(\tau),\ldots,\bm \varpi_{N}^a(\tau)]^T\\
\bm X^b(\tau)&=[\bm \gamma_1^b(\tau),\ldots,\bm \gamma_{N}^b(\tau),\bm \varpi_1^b(\tau),\ldots,\bm \varpi_{N}^b(\tau)]^T
\end{split}
\end{equation}
which are initialized by $\bm X^a(0)=\bm I_{2N}$ and $\bm X^b(0)=$ $[\bm \gamma_1^b(0),\ldots,\bm \gamma_{N}^b(0),\bm \varpi_1^b(0),\ldots,\bm \varpi_{N}^b(0)]^T$. The vectors $\bm \gamma_i^b(0)$ and $\bm \varpi_i^b(0)$ assigned to the $i^{th}$ bus satisfies:
\begin{small}
\begin{equation} \label{eq:vectorize}
\begin{split}
[{\bm \gamma_i^b(0)};{\bm \varpi_i^b(0)}]=[&{\mathbb{V}(\bm K_{\bm x_{Gi}}^{\bm f_{Gi}})};{\mathbb{V}(\bm K_{\bm \theta_{i}}^{\bm f_{Gi}})};{\mathbb{V}(\bm K_{\bm v_{i}}^{\bm f_{Gi}})};{\mathbb{V}(\bm K_{\bm x_{Gi}}^{\bm h_{p_{Gi}}})};\\
&{\mathbb{V}(\bm K_{\bm x_{Gi}}^{\bm h_{q_{Gi}}})};\bm 0],~i\in \mathcal{V}_G;\\
[{\bm \gamma_i^b(0)};{\bm \varpi_i^b(0)}]=[&{\mathbb{V}(\bm K_{\bm \theta_{i}}^{\bm f_{Wi}})};{\mathbb{V}(\bm K_{\bm x_{Wi}}^{\bm f_{Wi}})};{\mathbb{V}(\bm K_{\bm \omega_{i}}^{\bm f_{Wi}})};{\mathbb{V}(\bm K_{\bm v_{i}}^{\bm f_{Wi}})};\\
&{\mathbb{V}(\bm K_{\bm x_{Wi}}^{\bm h_{p_{Wi}}})};{\mathbb{V}(\bm K_{\bm \omega_{i}}^{\bm h_{p_{Wi}}})};{\mathbb{V}(\bm K_{\bm x_{Wi}}^{\bm h_{q_{Wi}}})};\\
&{\mathbb{V}(\bm K_{\bm \omega_{i}}^{\bm h_{q_{Wi}}})};\bm 0],~i\in \mathcal{V}_W;\\
[{\bm \gamma_i^b(0)};{\bm \varpi_i^b(0)}]=[&{\mathbb{V}(\bm K_{\bm \theta_i}^{\bm f_{Li}})};{\mathbb{V}(\bm K_{\bm x_{Li}}^{\bm f_{Li}})};{\mathbb{V}(\bm K_{\bm x_{Li}}^{\bm h_{p_{Li}}})};{\mathbb{V}(\bm K_{\bm \omega_{i}}^{\bm h_{p_{Li}}})};\\
&{\mathbb{V}(\bm K_{\bm \omega_{i}}^{\bm h_{q_{L_i}}})};\bm 0],~i\in \mathcal{V}_L;\\
[{\bm \gamma_i^b(0)};{\bm \varpi_i^b(0)}]=[&\bm 0],~i\in \mathcal{V}_T.
\end{split}
\end{equation}
\end{small}

The designed update process (\ref{eq:update}) can be realized in a distributed manner via the communication network $\mathcal{G}_1=\{\mathcal{V},\mathcal{E}\}$ mentioned earlier since
\begin{enumerate}
\item the initial values of vectors $\bm \gamma_i^a(0)$, $\bm \varpi_i^a(0)$, $\bm \gamma_i^b(0)$, and $\bm \varpi_i^b(0)$ can be assigned locally for each bus $i$, because i) the vectors $\bm \gamma_i^a(0)$, $\bm \varpi_i^a(0)$ can be assigned locally as each bus knows its own bus number and ii) the elements of vectors $\bm \gamma_i^b(0)$, $\bm \varpi_i^b(0)$ can be calculated based on local measurements $\theta_i$, $v_i$, $p_i^{inj}$ and $q_i^{inj}$ \cite{stability_control,wind_turbine};
\item for each sub-matrix $\bm J_{\bm \theta}^{\bm h_p}$, $\bm J_{\bm v}^{\bm h_p}$, $\bm J_{\bm \theta}^{\bm h_q}$, $\bm J_{\bm v}^{\bm h_q}$ of $\bm J_{pf}$ (see (\ref{eq:TMatrix}) and (\ref{eq:Spf}) in Appendix C for details), the nonzero elements of the $i^{th}$ row are functions of measurements of bus $i$ and its neighboring bus $j\in \mathcal{N}_i$ \cite{fully_distributed,fully_micro}.
\end{enumerate}

During the update process, at each step $\tau,~0<\tau\leq2N$, each damping bus $i,~i\in\mathcal{V}_W\cup\mathcal{V}_L$ stores its own data and data from its neighboring buses (which can be realized via communication links between neighboring buses). Thus, the whole distributed algorithm is expressed as:
\begin{equation} \label{eq:whole}
\begin{split}
&[\bm X^a(\tau+1),\bm X^b(\tau+1)]=\bm J_{pf} [\bm X^a(\tau),\bm X^b(\tau)]\\
&[\bm \xi_i^a(\tau),\bm \xi_i^b(\tau)]=\bm S_i[\bm X^a(\tau),\bm X^b(\tau)],~i\in \mathcal{V}_W\cup\mathcal{V}_L
\end{split}
\end{equation}
where the matrix $\bm S_i=[\bm e_i,\bm e_j,\bm e_{N+i},\bm e_{N+j}]^T\in \mathbb{R}^{2(\mathcal{D}_i+1)\times 2N}$, $j\in \mathcal{N}_i$ selects the rows with respect to the damping bus $i$ and its neighboring buses $j,~j\in \mathcal{N}_i$; $\bm \xi_i^a(\tau)$, $\bm \xi_i^b(\tau)\in \mathbb{R}^{2(\mathcal{D}_i+1)\times 2N}$ denote the data collected by the damping bus $i$. We assume the discrete-time system (\ref{eq:whole}) is observable, which usually holds in practice \cite{fully_micro}, i.e., $rank(\bm O_i)=2N$ where $\bm O_i\in \mathbb{R}^{4(\mathcal{D}_i+1)N\times 2N}$ is defined as
\begin{equation} \label{eq:O}
\bm O_i=[\bm S_i^T,{(\bm S_i\bm J_{pf})}^T,\ldots,{(\bm S_i\bm J_{pf}^{2N-1})}^T]^T.
\end{equation}

After the update process (\ref{eq:whole}), each damping bus $i,~i\in\mathcal{V}_W\cup\mathcal{V}_L$ can recover $\bm J_{pf}$ and $\bm X^b(0)$ via the data it collected $\bm \xi_i^a(\tau)$ and $\bm \xi_i^b(\tau),~\tau=0,1,\ldots,2N$. For simplicity, we define the following data matrices:
\begin{equation} \label{eq:data}
\begin{split}
&\bm \Xi_{i1}^a=[\bm \xi_i^a(0)^T,\ldots,\bm \xi_i^a(2N-1)^T]^T\in\mathbb{R}^{4(\mathcal{D}_i+1)N\times 2N}\\
&\bm \Xi_{i2}^a=[\bm \xi_i^a(1)^T,\ldots,\bm \xi_i^a(2N)^T]^T\in\mathbb{R}^{4(\mathcal{D}_i+1)N\times 2N}\\
&\bm \Xi_{i}^a=[{\bm \Xi_{i1}^a}^T,{\bm \Xi_{i2}^a}^T]^T\in\mathbb{R}^{8(\mathcal{D}_i+1)N\times 2N}\\
&\bm \Xi_{i1}^b=[\bm \xi_i^b(0)^T,\ldots,\bm \xi_i^b(2N-1)^T]^T\in\mathbb{R}^{4(\mathcal{D}_i+1)N\times 2N}.
\end{split}
\end{equation}
The singular value decomposition of $\bm \Xi_i^a$ is also needed, which is given as:
\begin{equation} \label{eq:de}
\bm \Xi_i^a=[\tilde{\bm U}_{\xi_i},\tilde{\bm U}_{\xi_i}^0]\left[
          \begin{array}{c}
          \bm \Sigma_{\xi_i}\\
          \bm 0
          \end{array}
        \right]\tilde{\bm V}_{\xi_i}^T=\tilde{\bm U}_{\xi_i}\bm \Sigma_{\xi_i}\tilde{\bm V}_{\xi_i}^T
\end{equation}
where $\bm \Sigma_{\xi_i}$, $\tilde{\bm V}_{\xi_i}\in \mathbb{R}^{2N\times2N}$, $\tilde{\bm U}_{\xi_i}\in \mathbb{R}^{8(\mathcal{D}_i+1)N\times 2N}$, $\tilde{\bm U}_{\xi_i}^0\in \mathbb{R}^{8(\mathcal{D}_i+1)N\times (8(\mathcal{D}_i+1)N-2N)}$. Based on the matrices given in (\ref{eq:data}) and (\ref{eq:de}), each damping bus $i,~i\in\mathcal{V}_W\cup\mathcal{V}_L$ can recover $\bm J_{pf}$ and $\bm X^b(0)$ by the following equations:
\begin{subequations}
\begin{equation} \label{eq:recovery1}
\bm J_{pf}=(\tilde{\bm U}_{\xi_{i1}}^T\bm \Xi_{i1}^a)^{-1}\Theta_i\tilde{\bm U}_{\xi_{i1}}^T\bm \Xi_{i1}^a
\end{equation}
\begin{equation} \label{eq:recovery2}
\bm X^b(0)=(\bm \Xi_{i1}^a)^{\dagger}\bm \Xi_{i1}^b
\end{equation}
\end{subequations}
where $\tilde{\bm U}_{\xi_{i1}}$, $\tilde{\bm U}_{\xi_{i2}}\in \mathbb{R}^{4(\mathcal{D}_i+1)N\times 2N}$ are sub-matrices of $\tilde{\bm U}_{\xi_i}$ with $\tilde{\bm U}_{\xi_i}=[\tilde{\bm U}_{\xi_{i1}}^T,\tilde{\bm U}_{\xi_{i2}}^T]^T$, $\Theta_i=(\tilde{\bm U}_{\xi_{i1}}^T\tilde{\bm U}_{\xi_{i2}})(\tilde{\bm U}_{\xi_{i1}}^T\tilde{\bm U}_{\xi_{i1}})^{-1}\in\mathbb{R}^{2N\times2N}$, and the superscript $\dagger$ denotes the Moore-Penrose inverse. The mathematical proof of (\ref{eq:recovery1})-(\ref{eq:recovery2}) can be found in \cite{fully_micro}.

As mentioned earlier, each damping bus is assumed to know the model structures of SG, WTG, and load. Thus, each damping bus can identify the type of bus $i$ (i.e, SG, WTG, load, or transfer bus) based on the $\bm \gamma_i^b(0)$ and $\varpi_i^b(0)$ of $\bm X^b(0)$ obtained, and hence can recover all $\bm K_\vee^\wedge$ Jacobian matrices in $\bm K_1$, $\bm K_2$, and $\bm K_3$ of (\ref{eq:reform}) from $\bm X^b(0)$ obtained based on (\ref{eq:vectorize}). Combined with $\bm J_{pf}$ obtained, each damping bus can reconstruct $\bm A_r$ by (\ref{eq:re}) and (\ref{eq:reform}). Therefore, the critical LFO can be calculated by applying eigenvalue analysis to $\bm A_r$ at each damping bus.

\begin{remark}
In the proposed update process (\ref{eq:whole}), we assume that the sum of the length of all vectorized $\bm K_{\vee}^{\wedge}$ matrices related to each type of bus (i.e., SG, WTG, load, or transfer bus) is less than the length of the data vectors $[{\bm \gamma_i^b};{\bm \varpi_i^b}],~i\in\mathcal{V}$ assigned for each type of bus that is $4N$ (refer to (\ref{eq:vectorize}) for details). If there exist one type of bus whose sum of the length of all vectorized $\bm K_{\vee}^{\wedge}$ matrices is more than $4N$, additional data vectors $\bm \gamma_i^c$, $\bm \varpi_i^c\in \mathbb{R}^{2N}$ are assigned for each bus to form the additional data matrix $\bm X^c\in \mathbb{R}^{2N\times 2N}$. For the type of bus whose sum of the length of all the vectorized $\bm K_{\vee}^{\wedge}$ matrices is more than $4N$, $[{\bm \gamma_i^c(0)};{\bm \varpi_i^c(0)}]$ is initialized by the remaining elements. For the other types of buses whose sum of the length of the vectorized $\bm K_{\vee}^{\wedge}$ matrices is less than $4N$, $[{\bm \gamma_i^c(0)};{\bm \varpi_i^c(0)}]$ is initialized by zeros. The additional data matrix $\bm X^c(0)$ can be recovered by each damping bus via the same way as the data matrix $\bm X^b(0)$ is recovered.
\end{remark}

\subsection{Controller parameters tuning module}

In order to guarantee an adequate stability margin, the damping ratio $\varsigma_c$ of the critical LFO $\lambda_c=\sigma_c+j\omega_c$ should satisfy $\varsigma_c\geq\varsigma^{\star}$ where $\varsigma_c=-\sigma_c/|\lambda_c|$ and $\varsigma^{\star}>0$ is the preset threshold. Once the damping ratio of the critical LFO is less than $\varsigma^{\star}$, the parameters of each DCU (i.e., $K_i$, $T_{1i}$, $T_{2i}$, $T_{3i}$, and $T_{4i}$) of each damping bus will be tuned coordinately to improve the damping ratio of the critical LFO.

\begin{figure}[!t]
\begin{center}
\includegraphics[width=2.3in,height=1.2in]{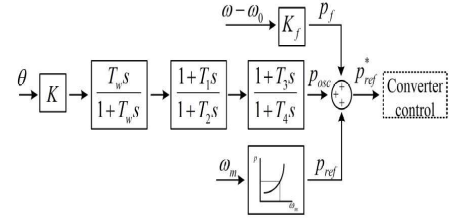}    
\caption{The closed-loop representation of the system.}
\label{fig:Hybrid}
\end{center}
\end{figure}

Without loss of generality, we firstly study the impact of the parameter changes of the $i^{th},~i\in\mathcal{V}_W\cup\mathcal{V}_L$ DCU (i.e., the DCU of the bus $N_G+i$) on $\lambda_c$. For analysis purposes, the system model (\ref{eq:ss}) is rewritten in the following form by reordering the variables of $\bm x$ in (\ref{eq:ss}):
\begin{equation} \label{eq:hybridss}
\left [\begin{array}{c}
\Delta \dot{\tilde{\bm x}}_i\\
\hdashline[2pt/2pt]
\bm 0
\end{array} \right]
=
\left [\begin{array}{c;{2pt/2pt}c}
\tilde{\bm A}_{si} & \tilde{\bm B}_{si}\\
\hdashline[2pt/2pt]
\tilde{\bm C}_{si} & \tilde{\bm D}_{si}\\
\end{array} \right]
\left [\begin{array}{c}
\Delta \tilde{\bm x}_i\\
\hdashline[2pt/2pt]
\Delta \bm y
\end{array} \right]
\end{equation}
where $\Delta \tilde{\bm x}_i=[\Delta \bm x_{i}^T,\Delta \bm x_{Ci}^T]^T$, $\bm x_{i}\in \mathbb{R}^{N_A-3}$ includes all state variables in $\bm x$ except $\bm x_{Ci}\in \mathbb{R}^{3}$ that is the corresponding state of the $i^{th}$ DCU; $\tilde{\bm A}_{si}=\bm T_i^{-1}\bm A_s\bm T_i$ (here $\bm T_i\in \mathbb{R}^{N_A\times N_A}$ is invertable which is the corresponding elementary row operator such that $\bm x=\bm T_i\tilde{\bm x}_i$); $\tilde{\bm B}_{si}=\bm T_i^{-1}\bm B_s$; $\tilde{\bm C}_{si}=\bm C_s\bm T_i$; and $\tilde{\bm D}_{si}=\bm D_s$. Then the system model (\ref{eq:hybridss}) can be written in the closed-loop form \cite{sensitivity}. In the closed-loop form, the system model is partitioned into two subsystems. For subsystem 1, which does not depend on parameters of the $i^{th}$ DCU, we have the following state space description:
\begin{equation} \label{eq:part10}
\left [\begin{array}{c}
\Delta \dot{\bm x}_{i}\\
\hdashline[2pt/2pt]
\bm 0
\end{array} \right]
=
\left [\begin{array}{c;{2pt/2pt}c}
\bm A_{i} & \bm B_{i}\\
\hdashline[2pt/2pt]
\bm C_{i} & \bm D_{i}\\
\end{array} \right]
\left [\begin{array}{c}
\Delta \bm x_{i}\\
\hdashline[2pt/2pt]
\Delta \bm y
\end{array} \right]
+
\left [\begin{array}{c}
\bm E_{i}\\
\hdashline[2pt/2pt]
\bm F_{i}
\end{array} \right]\Delta u_i.
\end{equation}
where $u_i=p_{osci}$ is the output of the $i^{th}$ DCU. Assuming $\bm D_i$ is nonsingular and eliminating $\Delta \bm y$ in (\ref{eq:part10}) gives:
\begin{equation}\label{eq:part1}
\Delta \dot{\bm x}_i=\bm A_{si}\Delta \bm x_i+\bm B_{si} \Delta u_i;~\Delta \theta_i=\bm C_{si}\Delta \bm x_i
\end{equation}
where $\bm A_{si}=\bm A_i-\bm B_i\bm D_i^{-1}\bm C_i\in \mathbb{R}^{(N_A-3)\times (N_A-3)}$, $\bm B_{si}=\bm E_i-\bm B_i\bm D_i^{-1}\bm F_i\in \mathbb{R}^{N_A-3}$ and $\bm C_{si}^T\in\mathbb{R}^{N_A-3}$. For subsystem 2, which only depends on the parameters of the $i^{th}$ DCU, we have the following state space description:
\begin{equation} \label{eq:part20}
\left [\begin{array}{c}
\Delta \dot{\bm x}_{Ci}\\
\hdashline[2pt/2pt]
\Delta u_i
\end{array} \right]
=
\left [\begin{array}{c;{2pt/2pt}c}
\bm A_{Ci} & \bm B_{Ci}\\
\hdashline[2pt/2pt]
\bm C_{Ci} & \bm D_{Ci}\\
\end{array} \right]
\left [\begin{array}{c}
\Delta \bm x_{Ci}\\
\hdashline[2pt/2pt]
\Delta \theta_i
\end{array} \right].
\end{equation}
where $\bm A_{Ci}$, $\bm B_{Ci}$, $\bm C_{Ci}$, and $\bm D_{Ci}$ can be easily obtained from (\ref{eq:DDCU}). A transfer function description for (\ref{eq:part20}) is given as:
\begin{equation} \label{eq:part2}
F_i(s,K_i)=\bm C_{Ci}(s\bm I-\bm A_{Ci})^{-1}\bm B_{Ci}+\bm D_{Ci}
\end{equation}
where $K_i\in \mathbb{R}$ is the gain factor in the $i^{th}$ DCU model. Based on (\ref{eq:part1}) and (\ref{eq:part2}), the schematic diagram of the closed-loop form is given in Fig. \ref{fig:Hybrid}.

Then the sensitivity of $\lambda_c$ with respect to $K_i$ of the transfer function $F_i(s,K_i)$ is given by \cite{sensitivity}:
\begin{equation} \label{eq:sen}
\frac{\partial \lambda_c}{\partial K_i}=\left.R_i\frac{\partial F_i(s,K_i)}{\partial K_i}\right|_{s=\lambda_c}
\end{equation}
where $R_i=\bm C_{si}\bm \phi_{si}\bm \psi_{si}^T \bm B_{si}\in \mathbb{C}$ is the residue with respect to the critical eigenvalue $\lambda_c$; $\bm \phi_{si}\in \mathbb{R}^{N_A-3}$ and $\bm \psi_{si}\in \mathbb{R}^{N_A-3}$ are the right and left eigenvectors of $\lambda_c$, respectively. Here, $\bm \phi_{si}$ ($\bm \psi_{si}$) consists of the first $N_A-3$ elements of $\bm \phi_i\in \mathbb{R}^{N_A}$ ($\bm \psi_i\in \mathbb{R}^{N_A}$) which is the right (left) eigenvector of $\lambda_c$ with respect to $\tilde{\bm A}_{ri}$ that is obtained by eliminating $\Delta\bm y$ in (\ref{eq:hybridss}), i.e.,
\begin{equation} \label{eq:hybridre}
\tilde{\bm A}_{ri}=\tilde{\bm A}_{si}-\tilde{\bm B}_{si}\tilde{\bm D}_{si}^{-1}\tilde{\bm C}_{si}=\bm T_i^{-1}\bm A_r\bm T_i.
\end{equation}

Combining (\ref{eq:sen}) and the transfer function of DCU given in Fig. \ref{fig:DDCU}, the sensitivity of $\lambda_c$ with respect to $K_i$ becomes:
\begin{equation} \label{eq:seni}
s_i=\frac{\partial \lambda_c}{\partial K_i}=R_i\cdot\frac{10\lambda_c}{1+10\lambda_c}\cdot\frac{1+T_{1i}\lambda_c}{1+T_{2i}\lambda_c}\cdot\frac{1+T_{3i}\lambda_c}{1+T_{4i}\lambda_c}.
\end{equation}
Here, the wash-out time constant of each DCU is assumed to be 10, i.e., $T_{wi}=10$.

It follows from (\ref{eq:seni}) that the tuning process of DCUs can be split into two parts: 1) tuning parameters of lead-lag compensation of the $i^{th}$ DCU such that $\angle s_i=180^\circ$; and then 2) tuning gain factors $K_i$ of all DCUs such that $\lambda_c$ moves to the desired location, i.e., $\sum_{i=1}^{N_W+N_L} |s_i|\Delta K_i\geq \Delta \Re(\lambda_c)^\star$ where $\Delta \Re(\lambda_c)^\star=\varsigma^\star|\omega_c|/\sqrt{1-(\varsigma^\star)^2}+\sigma_c$ is the expected real part change of $\lambda_c$. For part 1), the parameters of $T_{1i}$, $T_{2i}$, $T_{3i}$, and $T_{4i}$ can be calculated by \cite{sensitivity1}
\begin{equation} \label{eq:angle}
\begin{cases}
&\alpha_i=(1+\sin(\angle K_i)/2)/(1-\sin(\angle K_i)/2)\\
&T_{1i}=T_{3i}=(\sqrt{\alpha_i})/\omega_c\\
&T_{2i}=T_{4i}=1/(\sqrt{\alpha_i}\omega_c)
\end{cases}
\end{equation}
where $\angle K_i=180^\circ-\angle R_i$. For part 2), the gain factor change $\Delta K_i$ is calculated by solving the following optimization problem:
\begin{alignat}{2}\label{eq:obj}
\min~&\sum_{i=1}^{N_W+N_L} c_i\\\label{eq:constraint1}
\text{s.t.}~&\sum_{i=1}^{N_W+N_L} |s_i|\Delta K_i\geq \Delta \Re(\lambda_c)^\star\\\label{eq:constraint2}
&\Delta K_i^{min} \leq \Delta K_i\leq \Delta K_i^{max},i=1,\ldots,N_W+N_L
\end{alignat}
where $\Delta K_i^{min}$ and $\Delta K_i^{max}$ are the lower and upper bounds on the gain factor of the $i^{th}$ DCU, respectively. To account for the damping controller adjustments, in this work, we introduce a simple quadratic cost function for the $i^{th}$ damping bus which has been widely used in the literature (e.g., \cite{MPC}), i.e., $c_i=\pi_i\Delta K_i^2$ and $\pi_i>0$ is the cost parameter assigned for the $i^{th}$ damping bus. The objective (\ref{eq:obj}) is to minimize the total control cost. For convenience, the convex optimization problem (\ref{eq:obj})-(\ref{eq:constraint2}) is rewritten in a compact form as:
\begin{equation} \label{eq:compact_form}
\min_{\Delta \bm K} \sum_{i=1}^{N_W+N_L} c_i(\Delta K_i)~~~s.t.~~ g(\Delta \bm K)\leq 0,\Delta K_i\in \Delta \hat{K}_i
\end{equation}
where $\Delta \bm K=[\Delta K_1,\ldots,\Delta K_{N_W+N_L}]^T$ denotes the gain factor changes of $N_W+N_L$ DCUs; $g(\Delta \bm K)\leq 0$ represents the global constraint in (\ref{eq:constraint1}); $\Delta \hat{K}_i$ represents the local constraint in (\ref{eq:constraint2}).

As mentioned earlier, in this module, the proposed two-part tuning process will be realized in a distributed manner. For the first part tuning process, it is realized locally as the $R_i$ required of the $i^{th}$ damping bus can be obtained locally. It follows from (\ref{eq:sen}) that $R_i$ can be calculated by $\tilde{\bm A}_{ri}$, $\bm B_{si}$, and $\bm C_{si}$. For $\bm C_{si}$, based on (\ref{eq:part1}), it can be easily obtained as the $i^{th}$ damping bus knows the order of variables in $\bm x_i$. For $\tilde{\bm A}_{r_i}$, it can be calculated by (\ref{eq:hybridre}) as $\bm T_i$ is known locally and $\bm A_r$ has been obtained in the critical LFO identification module for each damping bus. For $\bm B_{si}$, based on (\ref{eq:hybridss})-(\ref{eq:part20}), we have
\begin{equation} \label{eq:relation}
\begin{split}
&\tilde{\bm A}_{si}
=
\left [\begin{array}{c;{2pt/2pt}c}
\bm A_{i}+\bm E_i\bm D_{Ci}\bm C_{si}   & \bm E_i\bm C_{Ci}\\
\hdashline[2pt/2pt]
\bm B_{Ci}\bm C_{si} & \bm A_{Ci}\\
\end{array} \right],~
\tilde{\bm B}_{si}
=
\left [\begin{array}{c}
\bm B_{i}\\
\hdashline[2pt/2pt]
\bm 0\\
\end{array} \right]\\
&\tilde{\bm C}_{si}
=
\left [\begin{array}{c;{2pt/2pt}c}
\bm C_{i}+\bm F_i\bm D_{Ci}\bm C_{si}   & \bm F_i\bm C_{Ci}
\end{array} \right],~
\tilde{\bm D}_{si}=\bm D_{i}.
\end{split}
\end{equation}
The $\bm B_{si}$ can be obtained by (\ref{eq:part1}) locally as: 1) matrices $\bm A_{Ci}$, $\bm B_{Ci}$, $\bm C_{Ci}$, and $\bm D_{Ci}$ is known locally, 2) according to (\ref{eq:hybridss}), matrices $\tilde{\bm A}_{si}$, $\tilde{\bm B}_{si}$, $\tilde{\bm C}_{si}$, and $\tilde{\bm D}_{si}$ can be calculated based on $\bm A_s$, $\bm B_s$, $\bm C_s$, and $\bm D_s$ which have been obtained by each damping bus in the critical LFO identification module, and 3) $\bm C_{si}$ can be obtained locally, then based on the matrix relations in (\ref{eq:relation}), matrices $\bm A_i$, $\bm B_i$, $\bm C_i$, $\bm D_i$, $\bm E_i$, and $\bm F_i$ can be calculated.

For the second part tuning process, in order to solve the convex optimization problem (\ref{eq:obj})-(\ref{eq:constraint2}) in a distributed manner, we decompose the Lagrange function of (\ref{eq:compact_form}) into a sum of $N_W+N_L$ local Lagrange functions where each of them is assigned to a damping bus:
\begin{equation} \label{eq:Lagrange}
L(\Delta \bm K,\mu)=\sum_{i=1}^{N_W+N_L} L_i(\Delta \bm K,\mu)
\end{equation}
where $L_i(\Delta \bm K,\mu)=c_i(\Delta K_i)+\mu g(\Delta \bm K)$, scalar $\mu$ is the Lagrange multiplier for $g(\Delta \bm K)\leq 0$ in (\ref{eq:compact_form}).

Inspired by (\ref{eq:Lagrange}), based on the distributed Lagrangian primal-dual sub-gradient algorithm proposed in \cite{DLPDS}, a distributed algorithm is designed to update the decision variables $\Delta \bm K$ and Lagrangian multiplier $\mu$ via communication between neighboring damping buses. The communication network used only covers damping buses and is allowed to have a different topology from the physical grid, which can be described by the undirected graph $\mathcal{G}_2=\{\mathcal{V}_2,\mathcal{E}_2,\mathcal{W}\}$, where $\mathcal{V}_2=\mathcal{V}_W\cup\mathcal{V}_L$, $\mathcal{E}_2\subseteq\mathcal{V}_2\times\mathcal{V}_2$, and $\mathcal{W}=\{w_{ij}\}\in\mathbb{R}^{(N_W+N_L)\times (N_W+N_L)}$. If $(i,j)\in \mathcal{E}_2$, $i \neq j$, then $w_{ij}=w_{ji}>0$ and $\sum_{j=1,i\neq j}^{N_W+N_L} w_{ij}<1$; otherwise, $w_{ij}=w_{ji}=0$. We define the diagonal entry $w_{ii}$ of the matrix $\mathcal{W}$ as $w_{ii}=1-\sum_{j=1,i\neq j}^{N_W+N_L} w_{ij}$. In the proposed distributed algorithm, the following assumptions are adopted:
\begin{enumerate}
\item The function $g$ in (\ref{eq:compact_form}) is known to all damping buses.
\item The topology of the communication network $\mathcal{G}_2$ is undirected and connected, and communication delays are negligible.
\end{enumerate}
For assumption 1), since $\bm A_s$, $\bm B_s$, $\bm C_s$, and $\bm D_s$ have been obtained by each damping bus via the critical LFO identification module, then all sensitivities $s_i$ in function $g$ can be calculated locally for each damping bus via the same method used for calculating $R_i$ in the first part tuning process.

Based on the abovementioned assumptions, the update process of decision variables $\Delta \bm K$ and Lagrangian multiplier $\mu$ is expressed as follows:
\begin{equation} \label{eq:update_dis}
\begin{split}
\Delta \bm K^{i}(\tau+1)&=P_{\Delta \hat K_i}[\Delta \bar{\bm K}^i(\tau)-\varsigma(\tau)\mathcal{D}L_{i,\Delta \bar{\bm K}^i}(\tau)]\\
\mu^i(\tau+1)&=P_{\hat U_i}[\bar{\mu}^i(\tau)+\varsigma(\tau)\mathcal{D}L_{i,\bar{\mu}^i}(\tau)]
\end{split}
\end{equation}
where $\Delta \bm K^{i}\in R^{N_W+N_L}$ and $\mu^i\in R$ are the information data assigned for the $i^{th}$ damping bus. We use $\Delta \bar{\bm K}^i(\tau)=\sum_{j=1}^{N_W+N_L} w_{ij}\Delta \bm K^j(r)$ and $\bar{\mu}^i(r)=\sum_{j=1}^{N_W+N_L} w_{ij}\mu^j(r)$ for short. At each time $\tau+1$, the $i^{th}$ damping bus calculates vectors $\mathcal{D}L_{i,\Delta \bar{\bm K}^i}=\partial L_i/\partial(\Delta \bar {\bm K}^i)$ and $\mathcal{D}L_{i,\bar{\mu}^i}=\partial L_i/\partial\bar{\mu}^i$ in the gradient direction of its local $L_i$. Combined with information received from its neighboring buses $\Delta \bar{\bm K}^i(\tau)$ and $\bar{\mu}^i(\tau)$, the $i^{th}$ damping bus updates its own decision variables $\Delta \bm K^i(\tau+1)$ and $\mu^i(\tau+1)$ by taking a projection onto its local constraint $\Delta \hat{K}_i$ and $\hat U_i=\{\mu_i\geq0\}$, respectively. Here, the projection operator $P_{\Delta \hat K_i}$ is defined by the definition of $P_{\Delta \hat K_i}[\bar{\bm x}]=\arg{\min}_{\bm x \in {\Delta \hat K_i}}\| \bar{\bm x}-\bm x\|$, where $\bar {\bm x}$ is a given vector. The projection operator $P_{\hat U_i}$ is defined in the same way as $P_{\Delta \hat K_i}$. The diminishing step size is $\varsigma(r)$ which satisfies $\lim_{r\to +\infty}\varsigma(r)=0$, $\sum_{r=0}^{+\infty}\varsigma(r)=+\infty$, and $\sum_{r=0}^{+\infty}\varsigma(r)^2<+\infty$. It has been proven in \cite{DLPDS} that for a convex optimization problem, the proposed distributed algorithm will asymptotically converge to a pair of primal-dual optimal solutions (i.e., $\lim_{\tau\to\infty} \Delta \bm K^i(\tau)=\Delta \bm K^\ast,~i=1,\ldots,N_W+N_L$ where $\Delta \bm K^\ast=[\Delta K_1^\ast,\cdots,\Delta K_{N_W+N_L}^\ast]^T$ is the optimal solution) under the Slater's condition, assumptions 1) and 2) mentioned above. In our case, the optimization problem (\ref{eq:obj})-(\ref{eq:constraint2}) is a convex optimization program whose global optimal solutions can be solved in a distributed way via the algorithm (\ref{eq:update_dis}). 

It is worth mentioning that, different damping buses have different geometric controllability/obserbility measures (COs) of the critical LFO $\lambda_c$ under different operating conditions \cite{H2}. The definition of the CO of the $i^{th}$ damping bus is given as $CO_i=\frac{|\bm \psi_{si}^T \bm B_{si}|}{||\bm \psi_{si}||||\bm B_{si}||}\cdot\frac{|\bm C_{si}\bm \phi_{si}|}{||\bm C_{si}||||\bm \phi_{si}||}$ which can be calculated locally. In the proposed two-part tuning process, only the damping buses with high COs participate the tuning process. In other words, if the CO of the $i^{th}$ damping bus satisfies $CO_i<CO^\star$ where $CO^\star$ is a threshold, then this damping bus does not participate the first part tuning process and the second part tuning process by setting $\Delta K_i^{min}=\Delta K_i^{max}=0$ in (\ref{eq:constraint2}) locally.

\section{Case study}

In this section, the modified IEEE 39-bus test system used for simulation is introduced firstly. Then, the simulation results and explanations will be presented.

\subsection{Test system}

\begin{figure}[!t]
\begin{center}
\includegraphics[width=2.6in,height=2.1in]{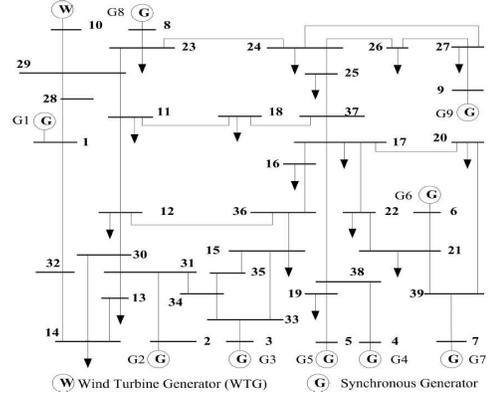}    
\caption{The modified IEEE 39-bus system}
\label{fig:IEEE39}
\end{center}
\end{figure}

Fig. \ref{fig:IEEE39} shows the modified IEEE 39-bus test system that is used to demonstrate the proposed distributed control framework. In the modified 39-bus system, the SG at bus 10 is replaced by a FRC-WTG with the same size of maximum power generation. All buses are renumbered according to the rules described in Section \uppercase\expandafter{\romannumeral2}-B, i.e., $N_G=9$, $N_W=1$, $N_L=17$, and $N_T=12$. The damping buses considered are 1 WTG bus and 17 load buses.
The model and system parameters are taken from \cite{ieee39}. For model parameters that are not provided in \cite{ieee39}, we use the default values of models given in the library developed in PSAT/MATLAB \cite{PSAT}.

For the communication network $\mathcal{G}_2$ used in the control parameters tuning module, each edge is assigned a weight which can be calculated by a simplified computational method:
\begin{equation} \label{eq:weight}
w_{ij}=\frac{1}{1+\max \{\tilde{\mathcal{D}}_i,\tilde{\mathcal{D}}_j\}},~i\in\mathcal{V}_2,~j\in \tilde{\mathcal{N}}_i
\end{equation}
where $\tilde{\mathcal{N}}_i$ defines the set of adjacent damping buses of the $i^{th}$ damping bus with the definition of $\tilde{\mathcal{N}}_i=\{j\in \mathcal{V}_2:(j,i)\in \mathcal{E}_2\}$ and cardinality $|\tilde{\mathcal{N}}_i|=\tilde{\mathcal{D}}_i$.

\subsection{Simulation results}

Following the procedure described in Section \uppercase\expandafter{\romannumeral3}-A, the critical LFO identified by each damping bus is $-0.0476\pm j1.7311$ with the damping ratio $\varsigma_c=0.0275$ (the preset threshold $\varsigma^{\star}=0.1$) and oscillation frequency equal to $0.275$ Hz, where G2-G9 oscillate against G1 (see Fig. \ref{fig:thetaG1}(a)). Then, the controller parameter tuning module is activated to tune the corresponding parameters as described in Section \uppercase\expandafter{\romannumeral3}-B. The price parameters needed for the optimization problem (\ref{eq:obj})-(\ref{eq:constraint2}) are given in Table \ref{tb:para_39}. In this case study, for simplicity, we assume the gain limits for DCUs are the same (i.e. $\Delta K^{min}=-10$ and $\Delta K^{max}=60$). As mentioned in Section \uppercase\expandafter{\romannumeral3}-B, only the damping buses with high COs participate the tuning process. Fig. 4 shows the COs of all 18 damping buses, and the threshold $CO^\star=10^{-4}$. It follows from Fig. 4 that buses 10, 17, 20, and 21 participate in the parameter tuning process. The obtained optimal gain factor changes are also given in Table \ref{tb:para_39}. The Table \ref{tb:lambda_c} compares the original $\lambda_c$, expected $\lambda_c$, and the new $\lambda_c$. It can be seen from Table \ref{tb:lambda_c} that the critical LFO is stabilized as desired.

\begin{figure}[!t]
\begin{center}
\includegraphics[width=3.0in,height=1.8in]{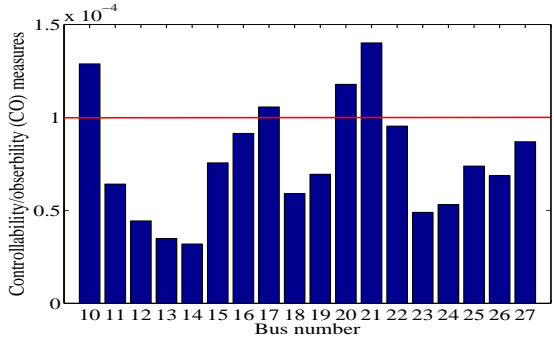}    
\caption{CO measures of all damping buses}
\label{fig:CO}
\end{center}
\end{figure}

\begin{table}[!t]
\begin{center}
\caption{The parameters for damping controllers}\label{tb:para_39}
\begin{tabular}{|c|c|c||c|c|c|}
\hline
Bus no.&$\pi$&$\Delta K^\star$&Bus no.&$\pi$&$\Delta K^\star$\\
\hline
10&0.8692&37&19&0.8524&0\\ \hline
11&0.9566&0&20&0.9367&60\\ \hline
12&0.7578&0&21&0.7306&33\\ \hline
13&1.2769&0&22&0.9391&0\\ \hline
14&0.8650&0&23&0.7993&0\\ \hline
15&1.3035&0&24&1.1363&0\\ \hline
16&1.3578&0&25&1.0443&0\\ \hline
17&0.9937&42&26&0.7898&0\\ \hline
18&1.0715&0&27&0.9862&0\\ \hline
\end{tabular}
\end{center}
\end{table}

\begin{table}[!t]
\begin{center}
\caption{Original, expected, and new $\lambda_c$}\label{tb:lambda_c}
\begin{tabular}{|c|c|c|}
\hline
Original $\lambda_c$&Expected $\lambda_c$&New $\lambda_c$\\ \hline
$-0.0476\pm j1.7311$&$-0.1749\pm j1.7311$&$-0.1785\pm j1.7340$\\ \hline
\end{tabular}
\end{center}
\end{table}

To illustrate the effectiveness of the proposed distributed control framework, we investigate the variation of rotor angle of G1 after a three-phase fault before and after the proposed tuning process. The three-phase fault happens at 1 s for 0.1 seconds on bus 25. From Fig. \ref{fig:thetaG1}(b) we can see that the system performance is improved with the presence of the proposed distributed control framework.

\begin{figure}[!t]
\begin{center}
\includegraphics[width=3.0in,height=2.0in]{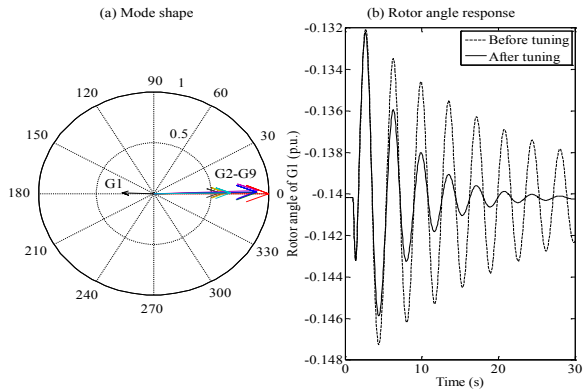}    
\caption{(a) Compass plot of relative mode shape (ref. G6) and (b) rotor angle of G1 responses to a three-phase fault on bus 25}
\label{fig:thetaG1}
\end{center}
\end{figure}

\section{Conclusion}

In this paper, WTGs and LAs have been coordinated to provide damping torques for the critical low frequency oscillation by adapting their active power generations and consumptions, respectively. In order to provide a scalable control framework for the increasing number of WTGs and LAs, a novel distributed control framework has been proposed which consists of a critical LFO identification module and a controller parameters tuning module. The simulation results have shown that the proposed distributed control framework is feasible and effective.

\section*{Appendix}

\subsection{SG model}

With $4^{th}$-order two-axis synchronous machine model and IEEE standard exciter model (IEEET1), the resulting differential-algebraic equations for the $i^{th}$ SG bus are given as:
\subsubsection{Differential equations}
\begin{equation} \label{eq:SG_d}
\begin{split}
\dot{e}'_{qi}=&\frac{1}{T'_{doi}}\left(-e'_{qi}-(x_{di}-x'_{di})i_{di}+v_{fi}\right)\\
\dot{e}'_{di}=&\frac{1}{T'_{qoi}}\left(-e'_{di}-(x_{qi}-x'_{qi})i_{qi}\right)\\
\dot{\delta}_i=&\omega_0(\omega_i-1)\\
\dot{\omega_i}=&\frac{1}{M_i}\left(p_{mi}-e'_{di}i_{di}-e'_{qi}i_{qi}-(x'_{qi}-x'_{di})i_{di}i_{qi}\right.\\
&\left.-D_i(\omega_i-1)\right)\\
\dot{v}_{mi}=&\frac{1}{T_{ri}}(v_i-v_{mi})\\
\dot{v}_{r1i}=&\frac{1}{T_{ai}}\left(K_{ai}(v_{refi}-v_{mi}-v_{r2i}-\frac{K_{fi}}{T_{fi}}v_{fi})-v_{r1i}\right)\\
\dot{v}_{r2i}=&-\frac{1}{T_{fi}}(\frac{K_{fi}}{T_{fi}}v_{fi}+v_{r2i})\\
\dot{v}_{fi}=&-\frac{1}{T_{ei}}\left(v_{fi}(K_{ei}+S_{ei}(v_{fi})-v_{ri})\right)
\end{split}
\end{equation}
where $\omega_0$ is the base frequency, $T'_{doi}$ and $T'_{qoi}$; $x_{di}$ and $x_{qi}$; $x'_{di}$ and $x'_{qi}$; $i_{di}$ and $i_{qi}$ are the d-axis and q-axis transient time constant; reactance; transient reactance; current, respectively; $p_{mi}$, $D_i$, and $M_i$ are the mechanical power, damping coefficient, and moment of inertia, respectively; $v_{fi}$ and $v_{ref_i}$ are the field and reference voltages, respectively; $T_{ri}$, $T_{ai}$, $T_{fi}$, and $T_{ei}$ are measurement, amplifier, stabilizer, and field circuit time constants, respectively; $K_{ai}$, $K_{fi}$, and $K_{ei}$ are amplifier, stabilizer, and field circuit gains, respectively; $S_{ei}$ is the ceiling function.
\subsubsection{Algebraic equations}
The stator algebraic equations are given as:
\begin{equation} \label{eq:SG_a1}
\begin{split}
p_{Gi}=&i_{di}v_i\sin(\delta_i-\theta_i)+i_{qi}v_i\cos(\delta_i-\theta_i)\\
q_{Gi}=&i_{di}v_i\cos(\delta_i-\theta_i)-i_{qi}v_i\sin(\delta_i-\theta_i).
\end{split}
\end{equation}
In order to express network voltages in the polar form, $i_{di}$ and $i_{q_i}$ in (\ref{eq:SG_d}) and (\ref{eq:SG_a1}) are expressed in terms of state variables $\bm x_{Gi}$ and algebraic variables $v_i$, $\theta_i$:
\begin{equation} \label{eq:SG_a2}
\left [\begin{array}{c}
i_{di}\\
i_{qi}
\end{array} \right]
=
\left [\begin{array}{cc}
r_{si}&-x'_{qi}\\
x'_{di}&r_{si}
\end{array} \right]^{-1}
\left [\begin{array}{c}
e'_{di}-v_i\sin(\delta_i-\theta_i)\\
e'_{qi}-v_i\cos(\delta_i-\theta_i)
\end{array} \right]
\end{equation}
where $r_{si}$ is the stator resistance. Substitution of (\ref{eq:SG_a2}) into (\ref{eq:SG_d}) and (\ref{eq:SG_a1}) gives \begin{equation}\label{eq:si_SI}
\begin{split}
\dot{\bm x}_{Gi}&=\bm f_{Gi}(\bm x_{Gi},\theta_i,v_i)\\
p_{Gi}&=g_{p_{Gi}}(\bm x_{Gi},\theta_i,v_i)\\
q_{Gi}&=g_{q_{Gi}}(\bm x_{Gi},\theta_i,v_i),~i\in\mathcal{V}_G
\end{split}
\end{equation}

\subsection{WTG model}

The model of a WTG includes models of the direct drive synchronous generator (DDSG), controller, and converter.
\subsubsection{DDSG model}
As the stator and rotor flux dynamics are fast in comparison with grid dynamics and the converter controls decoupled the generator from the grid, the steady-state electrical equations of DDSG are assumed. The differential and algebraic equations for DDSG of the $i^{th}$ WTG are given as:
\begin{equation} \label{eq:WTG_DDSG1}
\begin{split}
\dot{\omega}_{mi}&=\frac{1}{2H_{mi}}(\tau_{mi}-\tau_{ei})\\
p_{si}&=v_{sdi}i_{sdi}+v_{sqi}i_{sqi}\\
q_{si}&=v_{sqi}i_{sdi}-v_{sdi}i_{sqi}
\end{split}
\end{equation}
with
\begin{equation} \label{eq:WTG_DDSG2}
\begin{split}
\tau_{mi}&=\frac{p_{wi}(\theta_{pi})}{\omega_{mi}}\\
\tau_{ei}&=\psi_{sdi}i_{sqi}-\psi_{sqi}i_{sdi}\\
v_{sdi}&=-r_{si}i_{sdi}-\omega_{mi}\psi_{sqi}\\
v_{sqi}&=-r_{si}i_{sqi}+\omega_{mi}\psi_{sdi}\\
\psi_{sdi}&=-x_{sdi}i_{sdi}+\psi_{pmi}\\
\psi_{sqi}&=-x_{sqi}i_{sqi}
\end{split}
\end{equation}
where $H_{mi}$ is the rotor inertia; $p_{wi}(\theta_{pi})$ is the mechanical power which is the function of pitch angle $\theta_{pi}$; $\tau_{mi}$ and $\tau_{ei}$ are the mechanical and electrical torques, respectively; $v_{sdi}$ and $v_{sqi}$; $i_{sdi}$ and $i_{sqi}$; $x_{sdi}$ and $x_{sqi}$; $\psi_{sdi}$ and $\psi_{sqi}$ are stator d-axis and q-axis voltages; currents; reactances; and fluxes, respectively; $p_{si}$ and $q_{si}$ are produced active and reactive power, respectively; $r_{si}$ is the stator resistance; $\psi_{pmi}$ is the permanent magnet flux of rotor. Assuming the power factor equal to 1 (permanent magnet rotor), the reactive power output of the DDSG equals zero, i.e., $q_{si}=0$. Substituting (\ref{eq:WTG_DDSG2}) into (\ref{eq:WTG_DDSG1}) and expressing $i_{sqi}$ with $i_{sdi}$ based on $q_{si}=0$ in (\ref{eq:WTG_DDSG1}) gives:
\begin{equation} \label{eq:WTG_DDSG}
\begin{split}
\dot{\omega}_{mi}&=f_{mi}(\omega_{mi},\theta_{pi},i_{sqi})\\
p_{si}&=g_{spi}(i_{sqi}).
\end{split}
\end{equation}

\subsubsection{Controller}
The model of the controller includes models of the pitch angle control unit, primary frequency control unit, and the DDCU. For pitch angle control unit, its dynamic is described by the differential equation:
\begin{equation} \label{eq:WTG_pitch}
\dot{\theta}_{pi}=\frac{1}{T_{pi}}(K_{pi}\phi_i(\omega_{mi}-\omega_{mrefi})-\theta_{pi})
\end{equation}
where $K_{pi}$, $\omega_{mrefi}$, and $T_{pi}$ are pitch control gain, reference rotor speed, and pitch control time constant, respectively; $\phi_i$ is a function which allows varying the pitch angle set point only when the difference $\omega_{mi}-\omega_{mrefi}$ exceeds a predefined value $\pm\Delta\omega_{mi}$. For the primary frequency control unit, its control is given as:
\begin{equation} \label{eq:WTG_f}
p_{fi}=K_{fi}(\omega_i-\omega_0)
\end{equation}
with $\omega_i=\dot{\theta}_i$ where $\omega_0$ is the nominal frequency and $K_{fi}$ is the control gain. For the DDCU, its model is already given in Section \uppercase\expandafter{\romannumeral2}-A and repeated here for completeness:
\begin{equation}\label{eq:DDCU_SI}
\begin{split}
\dot{x}_{1i}=&-\frac{1}{T_{wi}}\left(K_i\theta_i+x_{1i}\right)\\
\dot{x}_{2i}=&\frac{1}{T_{2i}}\left((1-\frac{T_{1i}}{T_{2i}})(K_i\theta_i+x_{1i})-x_{2i}\right)\\
\dot{x}_{3i}=&\frac{1}{T_{4i}}\left((1-\frac{T_{3i}}{T_{4i}})\left(x_{2i}+\left(\frac{T_{1i}}{T_{2i}}(K_i\theta_i+x_{1i})\right)\right)-x_{3i}\right)\\
p_{osci}=&x_{3i}+\frac{T_{3i}}{T_{4i}}\left(x_{2i}+\frac{T_{1i}}{T_{2i}}(K_i\theta_i+x_{1i})\right).
\end{split}
\end{equation}
\subsubsection{Converter model}
Converter dynamics are highly simplified as they are fast in comparison with the electromechanical transients. Thus, the converter are modeled as an ideal current source where $i_{sq_i}$ and $i_{dc_i}$ are state variables and are used for the active power/speed control and the reactive power/voltage control, respectively. The differential equations for the converter of the $i^{th}$ WTG are given as:
\begin{equation} \label{eq:WTG_C1}
\begin{split}
\dot{i}_{sqi}&=\frac{1}{T_{pri}}(i_{sqrefi}-i_{sqi})\\
\dot{i}_{dci}&=\frac{1}{T_{Vi}}\left((v_{refi}-v_{i})-i_{cdi}\right)
\end{split}
\end{equation}
with
\begin{equation} \label{eq:WTG_C2}
i_{sqrefi}=\frac{p_{refi}(\omega_{mi})+p_{fi}+p_{osci}}{\omega_{mi}(\psi_{pmi}-x_{sdi}i_{sdi})}
\end{equation}
where $i_{sqrefi}$ is the reference current, $p_{refi}(\omega_{mi})$ is the power-speed characteristic which roughly optimizes the wind energy capture and is calculated by based on current rotor speed $\omega_{mi}$. The active and reactive power injected into the grid from the converter are given as:
\begin{equation} \label{eq:WTG_C3}
\begin{split}
p_{ci}&=v_{cdi}i_{cdi}+v_{cqi}i_{cqi}\\
q_{ci}&=v_{cqi}i_{cdi}-v_{cdi}i_{cqi}
\end{split}
\end{equation}
where $v_{cdi}=-v_{i}\sin\theta_i$ and $v_{cqi}=v_{i}\cos\theta_i$.

Assuming a lossless converter, the outputs of the WTG become
\begin{equation} \label{eq:WTG_C4}
\begin{split}
p_{Wi}&=p_{ci}=p_{si}\\
q_{Wi}&=v_i\left(i_{cdi}\cos\theta_i+\frac{\sin\theta_i(p_{si}+v_ii_{cdi}\sin\theta_i)}{v_i\cos\theta_i}\right).
\end{split}
\end{equation}
Substituting $p_{fi}$ in (\ref{eq:WTG_f}) and $p_{osci}$ in (\ref{eq:DDCU_SI}) into (\ref{eq:WTG_C2}), combining (\ref{eq:WTG_DDSG}), (\ref{eq:WTG_pitch}), (\ref{eq:DDCU_SI}), (\ref{eq:WTG_C1}), and (\ref{eq:WTG_C4}) gives
\begin{equation} \label{eq:ci_SI}
\begin{split}
\dot{\theta}_i&=\omega_i\\
\dot{\bm x}_{Wi}&=\bm f_{Wi}(\bm x_{Wi},\omega_i,\theta_i,v_i)\\
p_{Wi}&=g_{p_{Wi}}(\bm x_{Wi},\omega_i,\theta_i,v_i)\\
q_{Wi}&=g_{q_{Wi}}(\bm x_{Wi},\omega_i,\theta_i,v_i),~i\in\mathcal{V}_W.
\end{split}
\end{equation}

\subsection{Matrices}

\subsubsection{Matrices $\bm A_s$, $\bm B_s$, $\bm C_s$, and $\bm D_s$}

\begin{figure*}[!t]
\normalsize
\begin{multline}\label{eq:ABCD}
\left [\begin{array}{c;{2pt/2pt}c}
\bm A_s & \bm B_s\\
\hdashline[2pt/2pt]
\bm C_s & \bm D_s
\end{array} \right]=
\left[\arraycolsep=2pt{
               \begin{array}{ccccc;{2pt/2pt}cccccccc}
               \bm K_{\bm x_G}^{\bm f_G} & \bm 0&\bm 0&\bm 0&\bm 0&\bm K_{\bm \theta_G}^{\bm f_G} & \bm 0 & \bm 0 & \bm 0&\bm K_{\bm v_G}^{\bm f_G} & \bm 0 & \bm 0& \bm 0\\
               \bm 0 & \bm 0&\bm 0&\bm 0&\bm 0&\bm 0&\bm I_{N_W} & \bm 0 &\bm 0 & \bm 0&\bm 0 & \bm 0& \bm 0\\
               \bm 0 &\bm K_{\bm \theta_W}^{\bm f_W}& \bm K_{\bm x_W}^{\bm f_W}&\bm 0&\bm 0&\bm 0 & \bm K_{\bm \omega_W}^{\bm f_W} & \bm 0 & \bm 0 & \bm 0 &\bm K_{\bm v_W}^{\bm f_W} & \bm 0& \bm 0\\
               \bm 0 & \bm 0&\bm 0&\bm 0&\bm 0&\bm 0 &\bm 0 &\bm I_{N_L}& \bm 0&\bm 0 & \bm 0& \bm 0& \bm 0\\
               \bm 0 &\bm 0& \bm 0&\bm K_{\bm \theta_L}^{\bm f_L}&\bm K_{\bm x_L}^{\bm f_L}&\bm 0 & \bm 0 &\bm 0& \bm 0 & \bm 0& \bm 0& \bm 0& \bm 0\\\hdashline[2pt/2pt]
               \bm K_{\bm x_G}^{\bm h_{p_G}} & \bm J_{\bm \theta_W}^{\bm h_{p_G}} &\bm 0& \bm J_{\bm \theta_L}^{\bm h_{p_G}}&\bm 0&\bm J_{\bm \theta_G}^{\bm h_{p_G}}& \bm 0 & \bm 0 &\bm J_{\bm \theta_T}^{\bm h_{p_G}}& \bm J_{\bm v_G}^{\bm h_{p_G}} & \bm J_{\bm v_W}^{\bm h_{p_G}} & \bm J_{\bm v_L}^{\bm h_{p_G}}& \bm J_{\bm v_T}^{\bm h_{p_G}}\\
               \bm 0 & \bm J_{\bm \theta_W}^{\bm h_{p_W}} &\bm K_{\bm x_W}^{\bm h_{p_W}}& \bm J_{\bm \theta_L}^{\bm h_{p_W}}& \bm 0&\bm J_{\bm \theta_G}^{\bm h_{p_W}}&\bm K_{\bm \omega_W}^{\bm h_{p_W}} & \bm 0 &\bm J_{\bm \theta_T}^{\bm h_{p_W}}& \bm J_{\bm v_G}^{\bm h_{p_W}} & \bm J_{\bm v_W}^{\bm h_{p_W}} & \bm J_{\bm v_L}^{\bm h_{p_W}}& \bm J_{\bm v_T}^{\bm h_{p_W}}\\
               \bm 0 & \bm J_{\bm \theta_W}^{\bm h_{p_L}} &\bm 0& \bm J_{\bm \theta_L}^{\bm h_{p_L}}&\bm K_{\bm x_L}^{\bm h_{p_L}} &\bm J_{\bm \theta_G}^{\bm h_{p_L}}& \bm 0&\bm K_{\bm \omega_L}^{\bm h_{p_L}}&\bm J_{\bm \theta_T}^{\bm h_{p_L}}& \bm J_{\bm v_G}^{\bm h_{p_L}} & \bm J_{\bm v_W}^{\bm h_{p_L}} & \bm J_{\bm v_L}^{\bm h_{p_L}}& \bm J_{\bm v_T}^{\bm h_{p_L}}\\
               \bm K_{\bm x_G}^{\bm h_{q_G}} & \bm J_{\bm \theta_W}^{\bm h_{q_G}} &\bm 0& \bm J_{\bm \theta_L}^{\bm h_{q_G}}&\bm 0&\bm J_{\bm \theta_G}^{\bm h_{q_G}} & \bm 0 & \bm 0 &\bm J_{\bm \theta_T}^{\bm h_{q_G}}& \bm J_{\bm v_G}^{\bm h_{q_G}} & \bm J_{\bm v_W}^{\bm h_{q_G}} & \bm J_{\bm v_L}^{\bm h_{q_G}}& \bm J_{\bm v_T}^{\bm h_{q_G}}\\
               \bm 0 & \bm J_{\bm \theta_W}^{\bm h_{q_W}} &\bm K_{\bm x_W}^{\bm h_{q_W}}& \bm J_{\bm \theta_L}^{\bm h_{q_W}}&\bm 0&\bm J_{\bm \theta_G}^{\bm h_{q_W}}&\bm K_{\bm \omega_W}^{\bm h_{q_W}} & \bm 0 &\bm J_{\bm \theta_T}^{\bm h_{q_W}}& \bm J_{\bm v_G}^{\bm h_{q_W}} & \bm J_{\bm v_W}^{\bm h_{q_W}} & \bm J_{\bm v_L}^{\bm h_{q_W}}& \bm J_{\bm v_T}^{\bm h_{q_W}}\\
               \bm 0 & \bm J_{\bm \theta_W}^{\bm h_{q_W}} &\bm 0& \bm J_{\bm \theta_L}^{\bm h_{q_L}}&\bm 0&\bm J_{\bm \theta_G}^{\bm h_{q_L}} & \bm 0 &\bm K_{\bm \omega_L}^{\bm h_{q_L}} &\bm J_{\bm \theta_T}^{\bm h_{q_L}}& \bm J_{\bm v_G}^{\bm h_{q_L}} & \bm J_{\bm v_W}^{\bm h_{q_L}} & \bm J_{\bm v_L}^{\bm h_{q_L}}& \bm J_{\bm v_T}^{\bm h_{q_L}}
               \end{array}}
             \right]
\end{multline}
\begin{multline}\label{eq:TMatrix}
\left [\begin{array}{c;{2pt/2pt}c}
\bm K_1 & \bm K_2\\
\hdashline[2pt/2pt]
\bm K_3 & \bm J_{pf}
\end{array} \right]=
\left[\arraycolsep=2pt{
               \begin{array}{ccccc;{2pt/2pt}ccccccccc}
               \bm K_{\bm x_G}^{\bm f_G} & \bm 0&\bm 0&\bm 0&\bm 0&\bm K_{\bm \theta_G}^{\bm f_G} & \bm 0 & \bm 0 & \bm 0&&\bm K_{\bm v_G}^{\bm f_G} & \bm 0 & \bm 0&\bm 0\\
               \bm 0 & \bm I_{N_W}&\bm 0&\bm 0&\bm 0&\bm 0&\bm 0 & \bm 0 &\bm 0& &\bm 0 & \bm 0&\bm 0 &\bm 0\\
               \bm 0 &\bm K_{\bm \omega_W}^{\bm f_W}& \bm K_{\bm x_W}^{\bm f_W}&\bm 0&\bm 0&\bm 0 & \bm K_{\bm \theta_W}^{\bm f_W} & \bm 0 & \bm 0&&\bm 0 & \bm K_{\bm v_W}^{\bm f_W} & \bm 0& \bm 0\\
               \bm 0 & \bm 0&\bm 0&\bm I_{N_L}&\bm 0&\bm 0&\bm 0 &\bm 0 &\bm 0&&\bm 0 &\bm 0&\bm 0 & \bm 0\\
               \bm 0 &\bm 0& \bm 0&\bm 0&\bm K_{\bm x_L}^{\bm f_L}&\bm 0 & \bm 0 &\bm K_{\bm \theta_L}^{\bm f_L}&\bm 0 &&\bm 0 & \bm 0& \bm 0& \bm 0\\\hdashline[2pt/2pt]
               \bm K_{\bm x_G}^{\bm h_{p_G}} & \bm 0&\bm 0 & \bm 0 &\bm 0&\bm J_{\bm \theta_G}^{\bm h_{p_G}} & \bm J_{\bm \theta_W}^{\bm h_{p_G}}& \bm J_{\bm \theta_L}^{\bm h_{p_G}}&\bm J_{\bm \theta_T}^{\bm h_{p_G}} & \multicolumn{1}{r:}{}&\bm J_{\bm v_G}^{\bm h_{p_G}} & \bm J_{\bm v_W}^{\bm h_{p_G}} & \bm J_{\bm v_L}^{\bm h_{p_G}}& \bm J_{\bm v_T}^{\bm h_{p_G}}\\
               \bm 0 & \bm K_{\bm \omega_W}^{\bm h_{p_W}}&\bm K_{\bm x_W}^{\bm h_{p_W}} & \bm 0 &\bm 0&\bm J_{\bm \theta_G}^{\bm h_{p_W}} & \bm J_{\bm \theta_W}^{\bm h_{p_W}}& \bm J_{\bm \theta_L}^{\bm h_{p_W}} & \bm J_{\bm \theta_T}^{\bm h_{p_W}}& \multicolumn{1}{r:}{}&\bm J_{\bm v_G}^{\bm h_{p_W}} & \bm J_{\bm v_W}^{\bm h_{p_W}} & \bm J_{\bm v_L}^{\bm h_{p_W}}& \bm J_{\bm v_T}^{\bm h_{p_W}}\\
               \bm 0 & \bm 0&\bm 0 & \bm K_{\bm \omega_L}^{\bm h_{p_L}} &\bm K_{\bm x_L}^{\bm h_{p_L}}&\bm J_{\bm \theta_G}^{\bm h_{p_L}} & \bm J_{\bm \theta_W}^{\bm h_{p_L}}& \bm J_{\bm \theta_L}^{\bm h_{p_L}} & \bm J_{\bm \theta_T}^{\bm h_{p_L}}& \multicolumn{1}{r:}{}&\bm J_{\bm v_G}^{\bm h_{p_L}} & \bm J_{\bm v_W}^{\bm h_{p_L}} & \bm J_{\bm v_L}^{\bm h_{p_L}}& \bm J_{\bm v_L}^{\bm h_{p_L}}\\\cdashline{6-14}
               \bm K_{\bm x_G}^{\bm h_{q_G}} & \bm 0&\bm 0 & \bm 0 &\bm 0&\bm J_{\bm \theta_G}^{\bm h_{q_G}} & \bm J_{\bm \theta_W}^{\bm h_{q_G}}& \bm J_{\bm \theta_L}^{\bm h_{q_G}}&\bm J_{\bm \theta_T}^{\bm h_{q_G}}& \multicolumn{1}{r:}{}&\bm J_{\bm v_G}^{\bm h_{q_G}} & \bm J_{\bm v_W}^{\bm h_{q_G}} & \bm J_{\bm v_L}^{\bm h_{q_G}}& \bm J_{\bm v_T}^{\bm h_{q_G}}\\
               \bm 0 & \bm K_{\bm \omega_W}^{\bm h_{q_W}}&\bm K_{\bm x_W}^{\bm h_{q_W}} & \bm 0 &\bm 0&\bm J_{\bm \theta_G}^{\bm h_{q_W}} & \bm J_{\bm \theta_W}^{\bm h_{q_W}}& \bm J_{\bm \theta_L}^{\bm h_{q_W}}&\bm J_{\bm \theta_T}^{\bm h_{q_W}}& \multicolumn{1}{r:}{}&\bm J_{\bm v_G}^{\bm h_{q_W}} & \bm J_{\bm v_W}^{\bm h_{q_W}} & \bm J_{\bm v_L}^{\bm h_{q_W}}& \bm J_{\bm v_T}^{\bm h_{q_W}}\\
               \bm 0 & \bm 0&\bm 0 & \bm K_{\bm \omega_L}^{\bm h_{q_L}} &\bm 0&\bm J_{\bm \theta_G}^{\bm h_{q_L}} & \bm J_{\bm \theta_W}^{\bm h_{q_L}}& \bm J_{\bm \theta_L}^{\bm h_{q_L}}&\bm J_{\bm \theta_T}^{\bm h_{q_L}}& \multicolumn{1}{r:}{}&\bm J_{\bm v_G}^{\bm h_{q_L}} & \bm J_{\bm v_W}^{\bm h_{q_L}} & \bm J_{\bm v_L}^{\bm h_{q_L}}& \bm J_{\bm v_T}^{\bm h_{q_L}}\\\cdashline{6-14}
               \end{array}}
             \right]
\end{multline}
\end{figure*}

Refer to (\ref{eq:ABCD}) on next page for the detailed definition, where the notation $\bm K_{\vee}^{\wedge}$ ($\bm J_{\vee}^{\wedge}$) expresses Jacobian matrix of the $\vee$ in the subscript with respect to the $\wedge$ in the superscript. It should be noted that all the Jacobian matrices $\bm K_{\vee}^{\wedge}$ are block diagonal matrices.

\subsubsection{Matrices $\bm K_1$, $\bm K_2$, $\bm K_3$, and $\bm J_{pf}$}

Refer to (\ref{eq:TMatrix}) on next page for the detailed definition. All the Jacobian matrices $\bm J_\vee^\wedge$ in (\ref{eq:TMatrix}) form the power flow Jacobian matrix $\bm J_{pf}\in \mathbb{R}^{2N\times2N}$ where
\begin{equation}\label{eq:Spf}
\bm J_{pf}=\left[
          \begin{array}{c;{2pt/2pt}c}
           \bm J_{\bm \theta}^{\bm h_p} & \bm J_{\bm v}^{\bm h_p} \\
           \hdashline[2pt/2pt]
           \bm J_{\bm \theta}^{\bm h_q} & \bm J_{\bm v}^{\bm h_q}  \\
          \end{array}
        \right].
\end{equation}

\subsubsection{Elementary column operator matrix $\bm T$}
\begin{equation}\label{eq:T}
\left[\arraycolsep=0.1pt{\small{
               \begin{array}{ccccccccccccc}
               \bm I_{8N_G} & \bm 0&\bm 0 & \bm 0 & \bm 0 &\bm 0&\bm 0&\bm 0 & \bm 0 & \bm 0& \bm 0& \bm 0& \bm 0\\
               \bm 0 & \bm 0&\bm 0 & \bm 0 & \bm 0 &\bm 0&\bm I_{N_W}&\bm 0 & \bm 0 & \bm 0& \bm 0& \bm 0& \bm 0\\
               \bm 0 & \bm 0&\bm I_{7N_W} & \bm 0 & \bm 0 &\bm 0&\bm 0&\bm 0 & \bm 0 & \bm 0& \bm 0& \bm 0& \bm 0\\
               \bm 0 & \bm 0&\bm 0 & \bm 0 & \bm 0 &\bm 0&\bm 0&\bm I_{N_L} & \bm 0 & \bm 0& \bm 0& \bm 0& \bm 0\\
               \bm 0 & \bm 0&\bm 0 & \bm 0 & \bm I_{3N_L} &\bm 0&\bm 0&\bm 0 & \bm 0 & \bm 0& \bm 0& \bm 0& \bm 0\\
               \bm 0 & \bm 0&\bm 0 & \bm 0 & \bm 0 &\bm I_{N_G} &\bm 0&\bm 0 & \bm 0 & \bm 0& \bm 0& \bm 0& \bm 0\\
               \bm 0 & \bm I_{N_W}&\bm 0 & \bm 0 & \bm 0 &\bm 0 &\bm 0&\bm 0 & \bm 0 & \bm 0& \bm 0& \bm 0& \bm 0\\
               \bm 0 & \bm 0&\bm 0 & \bm I_{N_L} & \bm 0 &\bm 0 &\bm 0&\bm 0 & \bm 0 & \bm 0& \bm 0& \bm 0& \bm 0\\
               \bm 0 & \bm 0&\bm 0 & \bm 0 & \bm 0 &\bm 0 &\bm 0&\bm 0 & \bm I_{N_T} & \bm 0& \bm 0& \bm 0& \bm 0\\
               \bm 0 & \bm 0&\bm 0 & \bm 0 & \bm 0 &\bm 0 &\bm 0&\bm 0 & \bm 0 & \bm I_{N_G}& \bm 0& \bm 0& \bm 0\\
               \bm 0 & \bm 0&\bm 0 & \bm 0 & \bm 0 &\bm 0 &\bm 0&\bm 0 & \bm 0 & \bm 0& \bm I_{N_W}& \bm 0& \bm 0\\
               \bm 0 & \bm 0&\bm 0 & \bm 0 & \bm 0 &\bm 0 &\bm 0&\bm 0 & \bm 0 & \bm 0& \bm 0&\bm I_{N_L}& \bm 0\\
               \bm 0 & \bm 0&\bm 0 & \bm 0 & \bm 0 &\bm 0 &\bm 0&\bm 0 & \bm 0 & \bm 0& \bm 0 & \bm 0 &\bm I_{N_T}\\
               \end{array}}}
             \right]
\end{equation}

\end{document}